\documentclass{ifacconf}

\usepackage{graphicx}      
\usepackage{natbib}        
\usepackage{xcolor}
\usepackage[scr=boondox]{mathalpha}

\def\endpf{\hspace*{\fill}~\QED\par\endtrivlist\unskip}
\def\QEDopen{{\setlength{\fboxsep}{0pt}\setlength{\fboxrule}{0.2pt}\fbox{\rule[0pt]{0pt}{1.3ex}\rule[0pt]{1.3ex}{0pt}}}}
\def\QED{\QEDopen} 

\usepackage{epsfig} 
\usepackage{amsmath} 
\usepackage{amssymb}  

\begin{document}
\begin{frontmatter}

\title{{Indirect Adaptive Control Using a Static Update Law}
}
\author{Tom Kaufmann\thanksref{footnoteinfo1}}
\author{Johann Reger\thanksref{footnoteinfo2}} 
\thanks[footnoteinfo1]{
This author gratefully acknowledges support by Deutsche Forschungsgemeinschaft in the framework of Research Training Group “Tip- and laser-based $3$D-Nanofabrication in extended macroscopic working areas” at Technische Universität Ilmenau, Germany.
}
\thanks[footnoteinfo2]{Corresponding author}
\address{Control Engineering Group, Technische Universität Ilmenau, \mbox{P.O. Box $10$ $05$ $65$}, \mbox{D-$98684$ Ilmenau}, Germany,\\
e-mail: \{tom.kaufmann, johann.reger\}@tu-ilmenau.de}

\begin{abstract}
The update law in the indirect adaptive control scheme can be extended to include feedthrough of an error term. This reduces undesired oscillations of the calculated weights. When the $\sigma$-modification is used for achieving robustness against unstructured uncertainties, the gain of the feedthrough in the update law cannot be chosen arbitrarily. Compared to our previous result, we show stability of the closed loop for a larger parameter-range for the gain of the feedthrough in the update law. This parameter-range includes a configuration for which the influence of the integration in the update law diminishes over time, i.e. for which the adaptation for large times is governed solely by the feedthrough in the update law. By initializing at zero, this allows for removing the integration from the update law, resulting in a static update law. For the purely linear case, the adaptation acts like a disturbance observer. Frequency-domain analysis of the closed loop with a second order plant shows that removing the integration from the update law with $\sigma$-modification and feedthrough affects how precisely disturbances in the low-frequency band are observed. If the damping injected into the adaptation process by the $\sigma$-modification\linebreak exceeds certain bounds, then the precision is increased by using the static update law.
\end{abstract}
\end{frontmatter}
\vspace{-1.25mm}\section{Introduction}\vspace{-1.25mm}
Indirect adaptive control exploits known structures in the input disturbance of a nominally controlled error system to calculate weights with an associated update law. These weights are used in the adaptive control to counteract potentially destabilizing effects of the structured uncertainty on the nominal control loop. The classical approach based on the certainty equivalence principle in \cite{astrom_adaptive_2008} gives rise to an update law that integrates over an error term, which leads to drifting of the calculated weights when additional unstructured uncertainties are present in the input signal, see \cite{lavretsky_robust_2012}. There are several results available in the literature that face this lack of robustness. The \mbox{$\sigma$-modification} in \cite{ioannou_instability_1984} and \mbox{$\epsilon$-modification} in \cite{narendra_new_1987} both cause boundedness of the calculated weights by addition of different damping terms to the integrand of the update law. The projection method in \cite{pomet_adaptive_1989} proposes an extension of the classical approach that forces the calculated weights to remain in a predefined compact set.  According to \cite{weise_model_2024}, the performance of the closed loop when applying the \mbox{$\sigma$-modification} shows two major limitations: The error state generally does not converge to the origin even if no unstructured uncertainty is present, and high adaptation rates cause undesired oscillations of the calculated weights. Introducing feedthrough of an error term in the update law in \cite{weise_model_2024} helps reducing the oscillations without compromising the bandwidth of the adaptive control. However, the stability analysis in \cite{weise_model_2024} is carried out with regard to higher-order adaptive control, resulting in a conservative range for the gain of the feedthrough. \cite{kokcam_comparison_2024} prove zero residual error in the absence of unstructured uncertainty using an unmodified update law with feedthrough. We denote the update law with \mbox{$\sigma$-modification} and feedthrough in \cite{weise_model_2024} as proportional integral (PI) update law.

In this contribution, adaptation with a static update law is investigated as a special case of adaptation with the PI update law. After introducing both the error system and the sufficient assumptions in \mbox{Section $\ref{section:PROBLEM}$}, the PI update law is revisited in \mbox{Section \ref{part:PIUpdateLaw}} where stability is shown for a larger parameter-range. This allows the gain of the feedthrough to be chosen such that the integration can be removed from the PI update law, see \mbox{Section \ref{part:StaticUpdateLaw}}. Subsequently, some performance characteristics of the closed loop with the resulting static update law are investigated. \mbox{Section $\ref{part:ResidualError}$} highlights structured uncertainties which, provided there is no unstructured uncertainty, yield zero residual error. In \mbox{Section $\ref{part.GraphicalCriterion}$}, a criterion composed of conditions on the frequency-domain behavior of the error system is presented to ensure an increase in our estimate of the convergence rate of the error state when the gain of the static update law is increased. In \mbox{Section \ref{section:SIMULATION}}, the results are illustrated for a second order plant. First, the applicability of the criterion from \mbox{Section $\ref{part.GraphicalCriterion}$} is demonstrated. Second, we discuss the simulation results of the closed loop and, for the purely linear case, analyze the frequency-domain behavior of the adaptation. This shows that, despite its simpler structure, adaptation with the static update law achieves---without relying on high gain---a reduction of undesired oscillations of the calculated weights that is comparable to the effect of including the feedthrough in the PI update law. 

\textit{Notations}. The set of real symmetric matrices is denoted by \mbox{$\mathrm{S}^n=\left\{M\in\mathbb{R}^{n\times n}:\, M=M^\top\right\}$}. The matrix \mbox{$A^\mathrm{h}$} is the Hermitian of $A\in\mathbb{C}^{n\times m}$ and \mbox{$\|A\|_2=\sqrt{\lambda_{\max}\left(A^\mathrm{h}A\right)}$} is its maximal singular value, where $\lambda_{\min}(\cdot)$ and $\lambda_{\max}(\cdot)$ denote the minimal and maximal eigenvalue of any matrix in $\mathrm{S}^n$, respectively. The Euclidean norm of \mbox{$v\in\mathbb{C}^n$} is denoted by \mbox{$\|v\|=\|v\|_2$}. The identity of dimension $n$ is represented by $\mathrm{I}_n$. A matrix \mbox{$M\in\mathbb{R}^{n\times n}$} is called Hurwitz if all its eigenvalues are located in the open left complex half-plane. The Laplace transform of a time-domain signal is written as $\mathcal{L}(\cdot)$ dependent on $s\in\mathbb{C}$. We abbreviate “globally asymptotically stable” as GAS and “essentially bounded” as e.b. Whenever possible without causing misunderstandings, the arguments of a function are omitted for brevity.

\vspace{-1.25mm}\section{System and Assumptions} \label{section:PROBLEM} \vspace{-1.25mm}
We consider single-input error systems of the form
\begin{align}
	\label{ErrorDyn}
	\dot{e}=&Ae+B\, \big(u_\mathrm{adapt}+W^\top \beta(e)+\eta\big),
\end{align}
with known $A\!\in\!\mathbb{R}^{n\times n}$ Hurwitz, known $B\!\in\!\mathbb{R}^n$. The known function $\beta:\mathbb{R}^n\to\mathbb{R}^{n_\mathrm{\beta}}$ together with the unknown but constant weights $W\in\mathbb{R}^{n_\mathrm{\beta}}$ represents the structured uncertainty. Unstructured uncertainty is captured in $\eta(t)\in\mathbb{R}$. The uncertainties are matched to the input channel. We are mainly concerned with the following stability-concept:

\begin{defn}\cite{khalil_nonlinear_2002}
Let $x(t)\in\mathbb{R}^n$ with $x(t_0)=x_0$, $t_0\geq 0$ be the solution of a dynamical system that is defined in continuous time. The solution $x$ is called uniformly ultimately bounded (UUB) with the ultimate bound $b$ if $b>0$ exists independently of $t_0$ and $\forall a>0$: $\exists T=T(a,b)\geq 0$, independently of $t_0$, so that
	\begin{align}
		\|x(t_0)\|\leq a&\Longrightarrow\;\|x(t)\|\leq b\quad\forall t\geq t_0+T.
	\end{align}
\end{defn}
Further, let $P,Q\in\mathrm{S}^n$ satisfy the Lyapunov equation
\begin{align}
	\label{NomStab}
	Q=-\big(A^\top P+PA\big)
\end{align}
and, to enable the design of well-posed control, consider:
\begin{assum}\label{assumption:WellPosedness} There exists a non-decreasing function $\alpha_\mathrm{\beta}:[0,\infty)\to[0,\infty)\!$ with $\|\beta(x)\|\leq \alpha_\mathrm{\beta}\big(\|x\|\big)\!$ for all $x\in\mathbb{R}^n\!$.
\end{assum}
\begin{rem}
Nominal state feedback according to the model-reference-adaptive-control framework, as outlined in \cite{weise_model_2024}, yields the dynamics $(\ref{ErrorDyn})$ of the error between the plant and the reference model.\vspace{-0.75mm}
\end{rem}

\vspace{-1.25mm} \section{Results}\label{section:RESULTS} \vspace{-1.25mm}
The indirect adaptive control scheme aims at counteracting the influence of the structured uncertainty on the error system $(\ref{ErrorDyn})$ with the adaptive control effort
\begin{align}
    \label{ControlMRAC}
	u_\mathrm{adapt}=&-\hat{W}^\top\beta(e)
 \end{align}
 and the calculated weights $\hat{W}=\hat{W}(t)\in\mathbb{R}^{n_\beta}$, where we denote the algorithm to calculate $\hat{W}$ as update law.

\vspace{-1.25mm}\subsection{Proportional Integral Update Law}\label{part:PIUpdateLaw}\vspace{-1.25mm}
\begin{thm}\label{theorem:MRAC}
Consider the error dynamics $(\ref{ErrorDyn})$. Let Assumption $\ref{assumption:WellPosedness}$ apply and suppose, $0<P,Q\in\mathrm{S}^n$ satisfy the Lyapunov equation $(\ref{NomStab})$. The parameters $\Gamma,\Sigma,K\in\mathrm{S}^{n_\mathrm{\beta}}$ admit $0<\Gamma,\Sigma$ and $0\leq K<4\Sigma^{-1}$. Then, the adaptive control $(\ref{ControlMRAC})$ together with the PI update law\vspace{-1mm}
\begin{align}
	\label{UpdateLawMRAC}
	\hat{W}(t)\!=&K\beta(e(t))B^\top\!Pe(t)\!+\!\Gamma\!\!\int_0^t{\!\!\!\beta(e(\theta))B^\top\!Pe(\theta)\!-\!\Sigma \hat{W}(\theta)}\mathrm{d}\theta
\end{align}
guarantees UUB of $e$ and $\hat{W}$ if the unstructured uncertainty is bounded with $|\eta(t)|\leq \eta^\star<\infty$ for all $t\geq 0$.\vspace{-.75mm}
\end{thm}
\begin{pf}
Let
\begin{align}
\label{DefLyapunovFunction}
	\!\!\!\!V\!(e,\tilde{W})\!=\!e^\top\! Pe
	\!+\!(\tilde{W}\!-\!K\beta B^\top\! Pe)^{\!\top}\Gamma^{-1}\!(\tilde{W}\!-\!K\beta B^\top\! Pe)\!\!
\end{align}
with $\tilde{W}=\hat{W}-W\in\mathbb{R}^{n_\mathrm{\beta}}$ be a function. The closed-loop error dynamics read $\dot{e}=A e-B(\beta^\top\tilde{W})+B\eta$ and the update law can be expressed as its differentiated version $\frac{\mathrm{d}}{\mathrm{d}t}(\hat{W}-K\beta B^\top Pe)=\Gamma(\beta B^\top P e-\Sigma \hat{W})$. Therefore, the derivative of $V$ along the solution is
\begin{align}
\label{EqdotV}
	\!\dot{V}=&e^\top\big(A^\top P+PA-2PB\beta^\top K\beta B^\top P\big)e+2e^\top PB \eta\\\!\!&-2\tilde{W}^\top\Sigma\tilde{W}-2\tilde{W}^\top\Sigma {W}+2e^\top PB\beta^\top K\Sigma \hat{W}.\nonumber
\end{align}
For splitting off a square from the last term in $(\ref{EqdotV})$, we use that $2a^\top Kb\!=\!\delta a^\top K a\!+\!\delta^{-1}b^\top Kb\!-\!\delta^{-1}(b\!-\!\delta a)^\top K(b\!-\!\delta a)$ with $a\!=\!\beta B^\top Pe$ and $b\!=\!\Sigma \hat{W}$ for any scalar $\delta\!>\!0$. We get
\begin{align}
	\dot{V}=& e^\top\big(A^\top P+P A+(\delta-2)PB\beta^\top K\beta B^\top P\big)e\label{EqdotV2}\\
	&+2e^\top PB \eta+\tilde{W}^\top(\delta^{-1}\Sigma^\top K\Sigma-2\Sigma)\tilde{W}\nonumber\\
    &+2\tilde{W}^\top(\delta^{-1}\Sigma^\top K\Sigma-\Sigma)W+\delta^{-1} W^\top\Sigma^\top K\Sigma W\nonumber\\
	&-\delta^{-1}(\Sigma\hat{W}-\delta\beta B^\top P e)^\top K (\Sigma\hat{W}-\delta\beta B^\top P e)\nonumber
\end{align}
in view of $\hat{W}\!=\!\tilde{W}+W$. There exists $\delta^\star\!\in\!(0,2)$ with $2\Sigma^{-1}-{\delta^\star}^{-1}K\!>\!0$ since $2\Sigma^{-1}-2^{-1}K\!>\!0$. For $K\!\neq\! 0$, choose\linebreak e.g. $\delta^{\star-1}\!=\!2^{-1}\!+(2\lambda_{\mathrm{max}}(K))^{-1}\lambda_{\mathrm{min}}(2\Sigma^{-1}\!-2^{-1} K)$. Thus, $2\Sigma-\delta^{\star-1}\Sigma K\Sigma>0$ results from $0<\Sigma\in\mathrm{S}^{n_\mathrm{\beta}}$. We obtain
\begin{align}
    \!\dot{V}\!\leq\!& -c_0\|e\|^2\!-\!c_1\|\tilde{W}\|^2\!+\!2c_2\|e\|\!+\!2c_3\|\tilde{W}\|\!
	+\!\delta^{-1}W^\top\Sigma^\top\!K\Sigma W\nonumber\\	\!=\!&-c_0(\|e\|-c_0^{-1}c_2)^2-c_1(\|\tilde{W}\|-c_1^{-1}c_3)^2+c_4\label{IneqdotV}
\end{align}
when setting $\delta\!=\!\delta^\star$ and introducing
$c_0\!=\!\lambda_{\min}(Q)\!>\!0$, $c_1\!=\!\lambda_{\min}(2\Sigma\!-\!{\delta}^{-1}\Sigma
    K\Sigma)\!>\!0$, $c_3\!=\!\|(\delta^{-1}\Sigma^\top
    K\Sigma\!-\!\Sigma)W\|$, $c_2\!=\!\|PB\|\eta^\star$ and $c_4\!=\!c_0^{-1}c_2^2+c_1^{-1}c_3^2+\delta^{-1}W^\top \Sigma^\top K_\mathrm{p}\Sigma W$. Note the resulting implication
\begin{align}
            \label{InvariantRegion}
            \|e\|>v_1(\mu)\text{ or }\|\tilde{W}\|>v_2(\mu)&\Longrightarrow \dot{V}<-\mu
\end{align}
for $\mu\geq 0$ and functions $v_1(\mu)=c_0^{-1}c_2+(c_0^{-1}(c_4+\mu))^{1/2}$, $v_2(\mu)=c_1^{-1}c_3+(c_1^{-1}(c_4+\mu))^{1/2}$. Let the set
\begin{align}
\Omega_{\xi(\mu)}=&\big\{(e,\ \tilde{W})\in\mathbb{R}^{n}\times \mathbb{R}^{n_\beta}:\;V\leq \xi(\mu)\big\}
\end{align}
be depending on $\mu$ and function $\xi(\mu)=\lambda_{\max}(P)v_1(\mu)^2+\lambda_{\max}(\Gamma^{-1})(v_2(\mu)+\|K\|_2\alpha_\mathrm{\beta}(v_1(\mu))\|PB\|\, v_1(\mu))^2$.
The inversion of the implication $(\ref{InvariantRegion})$ yields
\begin{align}
	\dot{V}\geq 0&\Longrightarrow V\leq \xi(0)\leq\xi(\mu)\Longrightarrow (e,\;\tilde{W})\in\Omega_{\xi(\mu)},
\end{align}
which guarantees $\Omega_{\xi(\mu)}$ to be positively invariant. Provided there exists $T\geq 0$ with $V(e(T),\tilde{W}(T))\leq\xi(\mu)$, then $\Omega_{\xi(\mu)}$ being positively invariant implies \mbox{$\forall t\geq T: \|e(t)\|\leq r_1$} with the ultimate bound $r_1=(\xi(\mu)/\lambda_{\min}(P))^{1/2}$ of
$\|e\|$. Boundedness of $\hat{W}$ for all $t\geq T$ follows from the differentiated version of the PI update law: Addition of the term $0=-\Gamma\Sigma(-Kq)+\Gamma(-\Sigma K)q$ with $q=\beta B^\top Pe$ gives
\begin{align}
	\frac{\mathrm{d}}{\mathrm{d}t}(\hat{W}-Kq)=&-\Gamma\Sigma(\hat{W}-Kq)+\Gamma(\mathrm{I}_{n_\mathrm{\beta}}-\Sigma K)q \label{LTIVersionUpdateLaw}.
\end{align}
Equation~(\ref{LTIVersionUpdateLaw}) represents a linear time-invariant system with state $\hat{W}-Kq$ and input $q$. Its system matrix $-\Gamma\Sigma$ is Hurwitz since it shares all eigenvalues with $-\Gamma^{1/2}\Sigma\Gamma^{1/2}$, which in turn are all negative due to $0< \Gamma^{1/2}\in\mathrm{S}^n$ and $-\Sigma<0$. Thus, $\|e(T)\|\leq r_1$ implies for all $t\geq T$: The norm of the state $\|\hat{W}(t)-Kq(t)\|$ is bounded since the input admits $\|q(t)\|\leq \alpha_\mathrm{\beta}(r_1)\|PB\|r_1$. Therefore, $\|\hat{W}(t)\|$ is bounded as well. Now, we are only left to show the existence of $T\geq 0$ with $V(e(T),\tilde{W}(T))\leq\xi(\mu)$ in the case $V(e(0),\tilde{W}(0))>\xi(\mu)$. For a contradiction argument, suppose $V(e(t),\tilde{W}(t))>\xi(\mu)$ for all $t\geq 0$. Then, we have valid that $\|e(t)\|>v_1(\mu)$ or $\|\tilde{W}(t)\|>v_2(\mu)$ for each $t\geq 0$. Therefore, $\dot{V}(e(t),\tilde{W}(t))<-\mu$ for all $t\geq 0$. We obtain $V(e(t),\tilde{W}(t))< V(e(0),\tilde{W}(0))  -\mu t$ by means of the comparison principle. Fixing $\mu$ to a small positive value, this leads to the contradiction  $V(e(T),\tilde{W}(T))<\xi(\mu)$ for $T=(V(e(0),\tilde{W}(0))-\xi(\mu))\mu^{-1}>0$.\vspace{-.75mm}\end{pf}

\begin{rem}
    The conditions for UUB of $e$ and $\hat{W}$ in Theorem~\ref{theorem:MRAC} are less restrictive than the conditions required for the same result in \cite{weise_model_2024}, because here  neither $\Gamma,\Sigma$ diagonal nor $ K<\Sigma^{-1}$ is necessary.
\end{rem}

\vspace{-1.25mm}\subsection{Static Update Law}\label{part:StaticUpdateLaw}\vspace{-1.25mm}
Suppose, all conditions for UUB of $e$ and $\hat{W}$ in Theorem~\ref{theorem:MRAC} are met and choose\vspace{-1.5mm} 
\begin{align}
\label{FeedthroughGainForStaticUpdateLaw}
	K=\Sigma^{-1}
\end{align}  
as positive definite gain of the feedthrough 
in the PI update law $(\ref{UpdateLawMRAC})$. From $(\ref{EqdotV2})$ with $\delta=1$ we get
\begin{align}
	\label{easyDerivativeV}
\dot{V}\big|_{(\ref{FeedthroughGainForStaticUpdateLaw})}=&e^\top \big(A^\top P+PA-(\beta^\top K\beta) PBB^\top P\big)e\\
	&+2e^\top PB \eta-\tilde{W}^\top K^{-1}\tilde{W}+W^\top K^{-1} W\nonumber\\
	&-(\hat{W}-K\beta B^\top Pe)^\top K^{-1} (\hat{W}-K\beta B^\top Pe).\nonumber
\end{align}
Further, the differentiated PI update law from~(\ref{LTIVersionUpdateLaw}) turns into $\frac{\mathrm{d}}{\mathrm{d}t}(\hat{W}-K
    \beta B^\top Pe)\!=\!-\Gamma\Sigma(\hat{W}-K\beta B^\top
    Pe)$. Hence, $\lim_{t\to\infty}\hat{W}-K\beta B^\top Pe\!=\!0$ since $-\Gamma\Sigma$ is Hurwitz. Thereby:\vspace{-1mm}
\begin{align}
\label{VanishingIntegral}
\hat{W}=K\beta B^\top P e+ \underbrace{\Gamma \int_0^t \beta B^\top Pe-K^{-1} \hat{W}\mathrm{d}t}_{\to 0\text{ as } t\to\infty}.
\end{align}$\quad$\vspace{-1.4ex}\\
That is, the integral in the PI update law~(\ref{UpdateLawMRAC}) with~(\ref{FeedthroughGainForStaticUpdateLaw}) dies out after some time. Moreover, the static update law
\begin{align}
	\label{StaticUpdateLaw}
	\hat{W}(t)=K \beta(e(t))B^\top Pe(t)
\end{align}
holds for all times $t\geq 0$ if the integration in the PI update law
$(\ref{UpdateLawMRAC})$ with $(\ref{FeedthroughGainForStaticUpdateLaw})$ is initialized at zero. Note that
$(\ref{FeedthroughGainForStaticUpdateLaw})$ and $(\ref{StaticUpdateLaw})$ render the Lyapunov
function from $(\ref{DefLyapunovFunction})$
\begin{align}
	V\big|_{(\ref{StaticUpdateLaw})}=e^\top P e+W^\top \Gamma^{-1}W
\end{align}$\quad$\vspace{-2.25ex}\\
and its time derivative
\begin{align}
\label{EqdotVStaticUpdateLaw}
	\dot{V}\big|_{(\ref{StaticUpdateLaw})}=&e^\top \big(A^\top P+PA-(\beta^\top K\beta) PBB^\top P\big)e\\
	&+2e^\top PB \eta-\tilde{W}^\top K^{-1}\tilde{W}+W^\top K^{-1} W.\nonumber
\end{align}
\begin{rem}
    In the steady state, the PI update law~(\ref{UpdateLawMRAC}) with~(\ref{FeedthroughGainForStaticUpdateLaw}) behaves like the static update law~(\ref{StaticUpdateLaw}), compare~(\ref{VanishingIntegral}). Due to its lack of memory, adaptation with the static update law~(\ref{StaticUpdateLaw}) might be categorized as a \textit{robust} approach. However, we call it \textit{adaptive} since it represents a special case of adaptation with the PI update law~(\ref{UpdateLawMRAC}).\vspace{-0.75mm}
\end{rem}
In the following we investigate the behavior of the closed loop if the static
update law $(\ref{StaticUpdateLaw})$ is used to calculate the weights in the
adaptive control effort $(\ref{ControlMRAC})$. To this end, we give estimates of the ultimate bound, the residual error and the convergence rate of $\|e\|$, all in dependence of a linear scaling of the gain $K$ in $(\ref{StaticUpdateLaw})$, where
\begin{align}
\label{GainScaling}
	K=\alpha K_\mathrm{b}
\end{align}
with some constant base value $0\!<\!K_\mathrm{b}\!\in\!\mathrm{S}^{n_\mathrm{\beta}}$ and the scaling factor $\alpha>0$. In order to properly highlight the influence of the scaling $\alpha$ on those performance characteristics, for the remainder of Section~\ref{part:StaticUpdateLaw}, we assume the norm of $\beta$ to be lower bounded in the sense that $\mathscr{b}>0$ holds with 
\begin{align}
\label{LowerBoundBeta}
	\mathscr{b}=&\inf_{x\in\mathbb{R}^n}\big\{\beta(x)^\top K_\mathrm{b} \beta(x)\big\}.
\end{align}
$\quad$\vspace{-1.4ex}\\
Let $\mathrm{g}:[0,\infty)\times\mathrm{S}^{n}\times\mathbb{R}^{n\times m}\to\mathbb{R}$ be a function with
\begin{align}
g\big(\phi,Q,v\big)=&\lambda_{\min}\big(Q+\phi\; vv^\top\big)
\end{align}
and consider the virtual output $y_\mathrm{v}=B^\top P e$ as well as the scalar constant $\gamma\in[0,1)$. Completion of squares used on an upper bound of $(\ref{EqdotVStaticUpdateLaw})$ yields
\begin{align}
    &\dot{V}\big|_{(\ref{StaticUpdateLaw})}\leq-g\big(\alpha\mathscr{b}\gamma,Q,PB\big)\|e\|^2-\alpha\mathscr{b}(1-\gamma)\|y_\mathrm{v}\|^2\nonumber\\
    &\qquad\quad\;\;+2\|y_\mathrm{v}\|\eta^\star\!+\frac{1}{\alpha}W^\top K_\mathrm{b}^{-1} W\nonumber\\
    \label{IneqdotVStaticUpdateLaw}
    &\leq\! -\!\underbrace{g\big(\alpha\mathscr{b}\gamma,\!Q,\!PB\big)}_{\geq \lambda_{\min}(Q)>0}\|e\|^2\!+\!\frac{1}{\alpha}\!\left(\!W^\top \!K_\mathrm{b}^{-1}W\!+\!\frac{{\eta^\star}^2}{\mathscr{b}(1\!-\!\gamma)}\!\right)\!,\!\!
\end{align}
$\quad$\vspace{-1.4ex}\\
which describes the dependence of the closed-loop-behavior on $\alpha$ and, indirectly, on function $g$. In order to gain further insight, we list some properties of function $g$ in Lemma~\ref{lemma.EigenvalueLift}.

\begin{lem}\label{lemma.EigenvalueLift}
Let $Q\in\mathrm{S}^n$, $v\in\mathbb{R}^{n\times m}$ and $m<n$. Function $g(\cdot,Q,v)$ is non-decreasing, Lipschitz continuous with Lipschitz constant $\lambda_{\mathrm{max}}\big(vv^\top)$ and it remains in the bounds $\lambda_{\min}(Q)\leq g(\phi,Q,v)\leq\lambda_{\max}(Q)$ for all $\phi\geq 0$. Further:
\begin{align}
\label{EigenvalueLift}
\!\sup_{\phi\geq 0}\big\{g(\phi,Q,v)\big\}\!>\!\lambda_{\min}(Q)&\!\Longleftrightarrow\!\;g(1,Q,v)\!>\!\lambda_{\min}(Q). \!
\end{align}
$\quad$\vspace{-1.3ex}\\
Consider $\mathcal{N}=\mathrm{arg}\max_{\phi\geq 0}\big\{g(\phi,Q,v)\big\}$. Let $\phi^\star=\mathrm{inf}\mathcal{N}$ if $\mathcal{N}$ is nonempty, and $\phi^\star=\infty$ otherwise. Provided the statements in the equivalence $(\ref{EigenvalueLift})$ are true, then $\phi^\star>0$ and $g(\cdot,Q,v)\!$ is strictly increasing on $\mathcal{I}=[0,\phi^\star)$, i.e.
\begin{align}
\label{StrictlyIncreasing}
a,b\in \mathcal{I},\; a<b&\Longrightarrow g(a,Q,v)<g(b,Q,v),\\
\label{OnceNonStrictAlwaysNonStrict}
c\geq \phi^\star&\Longrightarrow g(c,Q,v)=\sup_{\phi\geq 0}\big\{g(\phi,Q,v)\big\}.
\end{align}
$\quad$\vspace{-1.3ex}\\Conversely, $g(\cdot,Q,v)\!\equiv\!\lambda_\mathrm{min}(Q)$ if $g(1,Q,v)\!=\!\lambda_\mathrm{\min}(Q)$.\vspace{-0.75mm}
\end{lem}

\begin{pf} Weyl's inequality, see Theorem~4.3.1 in \cite{horn_matrix_1994}, yields $\lambda_\mathrm{min}(Q)+\phi\,  \lambda_\mathrm{min}\big(vv^\top\big)\!\leq\!g(\phi,\! Q,\! v)\!\leq\! \lambda_{\mathrm{max}}(Q)+\phi\, \lambda_\mathrm{min}\big(vv^\top\big)$ for all $\phi\!\geq\! 0$. We are concerned with $m\!<\!n$, i.e. $\lambda_{\min}\big(vv^\top\big)\!=\!0$. Thus, $g(\cdot,\!Q,\!v)$ is bounded as claimed. Consider any $\phi,\chi\!\geq\! 0$ and let $w_1\!\in\!\mathbb{R}^n$ (or $w_2\!\in\!\mathbb{R}^n$) be a normalized eigenvector with respect to the minimal eigenvalue of $Q+(\phi+\chi) vv^\top$ (or $Q+\phi\, vv^\top$). The function $g(\cdot,\!Q,\!v)$ is non-decreasing since $g(\phi+\chi,\!Q,\!v)\!=\!w_1^\top\big(Q+\phi\, vv^\top\big)w_1+\chi\, \|v^\top\! w_1\|^2\!\geq\! g(\phi,\!Q,\!v)$. From $g(\phi,\!Q,\!v)\!=\!w_2^\top \big(Q+\phi\, vv^\top\big)w_2$ follows $g(\phi,\!Q,\!v)+\chi\, w_2^\top\! vv^\top\! w_2\!=\!w_2^\top\!\big(Q+(\phi+\chi)vv^\top\big)w_2\!\geq\! g(\phi+\chi,\!Q,\!v)$. Rearranging yields $g(\phi+\chi,\!Q,\!v)-g(\phi,\!Q,\!v)\!\leq\! \chi\, w_2^\top\! vv^\top\! w_2\!\leq\!\chi\, \lambda_{\mathrm{max}}\big(vv^\top)$. Since $g(\cdot,\!Q,\!v)$ is non-decreasing, its Lipschitz continuity with Lipschitz constant $\lambda_{\mathrm{max}}\big(vv^\top)$ results. The equivalence 
\begin{align}
\label{EigenvalueLiftProof}
g(\phi,Q,v)=\lambda_{\min}(Q)&\Longleftrightarrow\; g(1,Q,v)=\lambda_{\min}(Q)
\end{align}
for any $\phi\!>\!0$ follows from $g(1,\!Q,\!\sqrt{\phi}v)\!=\!g(\phi,\!Q,\!v)$~by~invoking  Lemma~\ref{lemma:MinimalEigenvalue} from the Appendix once with $w\!=\!\sqrt{\phi} v$ and once with $w\!=\!v$. Necessity in $(\ref{EigenvalueLift})$ is obvious. Sufficiency in\linebreak $(\ref{EigenvalueLift})$ arises via negation of the statements in~$(\ref{EigenvalueLiftProof})$,~since $g(\cdot,\!Q,\!v)$ is non-decreasing and, given $\sup_{\phi\geq 0}\{g(\phi,\!Q,\!v)\}\!>\!\lambda_{\min}(Q)$, there is $\phi'\!>\!0$ with~$g(\phi'\!,\!Q,\!v)\!>\!\lambda_{\min}(Q)$. Let the statements in~(\ref{EigenvalueLift}) be true. If $\phi^\star\!<\!\infty$, then $\mathcal{N}\!\neq\!\emptyset$~and, due to Lemma~\ref{lemma:WellPosednessLipschitzianFunction} from the Appendix, $\phi^\star\!\in\! \mathcal{N}$ since $g(\cdot,\!Q,\!v)$~is Lipschitz continuous. Thus, $\phi^\star\!>\!0$ because either $\phi^\star\!\!=\!\infty$ or $g(\phi^\star\!,\!Q,\!v)\!=\!\mathrm{sup}_{\phi\geq 0}\{g(\phi,\!Q,\!v)\}\!>\!g(0,\!Q,\!v)$. Implication~(\ref{OnceNonStrictAlwaysNonStrict}) results from $g(\cdot,\!Q,\!v)$ being non-decreasing and the construction of $\mathcal{I}$. Consider $\phi_2\!>\!\phi_1\!\geq\! 0$, $\chi\!>\!0$. Equivalence $(\ref{EigenvalueLiftProof})$ results in $g\big(\chi,\!Q+\phi_1\, vv^\top\!,\!v\big)\!=\!\lambda_{\min}\big(Q+\phi_1\, vv^\top\big)$ if and only if $g\big(1,\!Q+\phi_1\, vv^\top\!,\!v\big)\!=\!\lambda_{\min}\big(Q+\phi_1\, vv^\top\big)$, holding, if and only if $g\big(\phi_2-\phi_1,\!Q+\phi_1\, vv^\top\!,\!v\big)\!=\!\lambda_{\min}\big(Q+\phi_1\, vv^\top\big)$. Observing $g\big(\chi,\!Q+\phi_1\, vv^\top\!,\!v\big)\!=\!g(\phi_1+\chi,\!Q,\!v)$ and $g\big(\phi_2-\phi_1,\!Q+\phi_1\, vv^\top\!,\!v\big)\!=\!g(\phi_2,\!Q,\!v)$ leads to: If there is $\chi\!>\!0$ with $g(\phi_1+\chi,\!Q,\!v)\!=\!g(\phi_1,\!Q,\!v)$ for $\phi_1\!\geq\! 0$, then $g(\phi_2,\!Q,\!v)\!=\!g(\phi_1,\!Q,\!v)$ for all $\phi_2\!>\!\phi_1$. To prove $(\ref{StrictlyIncreasing})$ by contradiction, suppose $a,b\!\in\! \mathcal{I}$ with $a\!<\!b$ and $g(a,\!Q,\!v)\!=\!g(b,\!Q,\!v)$ exist. With $\phi_1\!=\!a$, $\chi\!=\!b-a$, we get $g(a^\star\!,\!Q,\!v)\!=\!g(a,\!Q,\!v)$ for all $a^\star\!\in\!\big\{\phi\!\in\! \mathcal{I}: \phi\!> \!a\big\}\!\neq\!\emptyset$ and thereby $\sup_{\phi> a}\{g(\phi,\!Q,\!v)\}\!=\!g(a,\!Q,\!v)$ by means of $(\ref{OnceNonStrictAlwaysNonStrict})$. Since $g(\cdot,\!Q,\!v)$ is non-decreasing, this implies \mbox{$\mathrm{sup}_{\phi\geq 0}\{g(\phi,\!Q,\!v)\}\!=\!g(a,\!Q,\!v)$ and thus $a\!\geq\!\phi^\star\!$, i.e. $a\!\notin\! \mathcal{I}$.}
\end{pf}

\vspace{-0.75mm}\subsubsection{Ultimate Bound}
For any scalar $\mu>0$ the number
\begin{align}   
\label{UltimateBound}
r_\mathrm{e}=&\sqrt{\frac{\lambda_{\max}(P)}{g\big(\alpha\mathscr{b}\gamma,Q,PB\big)\lambda_{\min}(P)}}\sqrt{\frac{W^\top K_\mathrm{b}^{-1}W+\frac{{\eta^\star}^2}{\mathscr{b}(1-\gamma)}}{\alpha}+\mu}
\end{align}
 is an ultimate bound of $\|e\|$ if the static update law is used in the closed loop to calculate $\hat{W}$. The constant $\mu$ determines how long it takes $\|e\|$ to reach $r_\mathrm{e}$ after its initial condition at $t=0$, with larger values for $\mu$ leading  to a shorter settling time.\vspace{-0.75mm}

\begin{pf}
Let $v_\mathrm{e}\!=\!\sqrt{\lambda_{\min}(P)/\lambda_{\max}(P)}\,r_\mathrm{e}$ depend on $\mu\!\geq\! 0$\linebreak and consider $\Omega^\mathrm{e}_{\rho(\mu)}\!=\!\big\{e\!\in\!\mathbb{R}^n:V\!\leq\! \rho(\mu)\big\}$
with $\rho(\mu)\!=\!\lambda_{\max}(P)v_\mathrm{e}(\mu)^2\!+\!W^\top\Gamma^{-1}W$. From $(\ref{IneqdotVStaticUpdateLaw})$, we get $\|e\|\!>\!v_\mathrm{e}(\mu)\Longrightarrow \dot{V}\big|_{(\ref{StaticUpdateLaw})}\!<\!-\mu$. Negating the statements of this implication yields $\dot{V}\!\geq\! 0\Longrightarrow \|e\|\!\leq\! v_\mathrm{e}(0)\!\leq\! v_\mathrm{e}(\mu)\Longrightarrow e\!\in\!\Omega_{\rho(\mu)}^\mathrm{e}$, which in turn exhibits the positive invariance of $\Omega_{\rho(\mu)}^\mathrm{e}$ for any $\mu\!\geq\! 0$. Further, if we restrict $\mu$ to positive values and consider $e(0)\!\notin\! \Omega_{\rho(\mu)}^\mathrm{e}$, we can guarantee $e(t)\!\in\!\Omega_{\rho(\mu)}^\mathrm{e}$ for all $t\!\geq\! T_\mathrm{e}$ with the settling time $T_\mathrm{e}\!=\!(V(e(0),\tilde{W}(0))|_{(\ref{StaticUpdateLaw})}\!-\!\rho(\mu))\mu^{-1}$ due to an argument similar to the contradiction used to prove Theorem~\ref{theorem:MRAC}. Observing $\|e\|\!\leq\! r_\mathrm{e}$ for all $e\in\Omega_{\rho(\mu)}^\mathrm{e}$ concludes this proof.\vspace{-1mm}
\end{pf}

The ultimate bound $r_\mathrm{e}\big|_{\mu>0}$ can be rendered arbitrarily small by selecting sufficiently large values for $\alpha$ and small values for $\mu>0$, i.e. by increasing the gain to the static update law accordingly and accepting long settling times.

\vspace{-0.75mm}\subsubsection{Transient Behavior}
If $\|e(0)\|\geq \left.r_\mathrm{e}\right|_{\mu=0}$, then the transient error is bounded as per
\begin{align}
\label{SlowestConvergence}
    \!\!\!\!\|e(t)\|\!\leq&\sqrt{r_\mathrm{e}\big|_{\mu=0}^2\!\!+\!\frac{1}{\mathrm{exp}(c_\mathrm{e} t)}\!\!\left(\!\frac{\lambda_{\max}(P)}{\lambda_{\min}(P)}\|e(0)\|^2\!-\!r_\mathrm{e}\big|_{\mu=0}^{2}\!\right)\!}\!\!
\end{align}
for all $t\geq 0$, with the convergence rate
\begin{align}
    c_\mathrm{e}=&g\big(\alpha\mathscr{b}\gamma, Q, PB\big)\big/\lambda_{\max}(P)
\end{align}
and the residual error $r_\mathrm{e}\big|_{\mu=0}$ of $\|e\|$.\vspace{-1mm}

\begin{pf} For this proof, the constant $\mu$ in $r_\mathrm{e}$ and $v_\mathrm{e}$ is set to\linebreak zero. We introduce the function  $U(e)=e^\top Pe-\lambda_{\min}(P){r_\mathrm{e}}^2$ and show that $\dot{U}\leq-c_\mathrm{e}U$ holds if $U\geq 0$. From $(\ref{IneqdotVStaticUpdateLaw})$ follows
{$\dot{U}=\dot{V}\big|_{(\ref{StaticUpdateLaw})}\leq -g_\mathrm{e}\, (\|e\|^2-{v_\mathrm{e}}^2)$} with $g_\mathrm{e}\!=\!g\big(\alpha\mathscr{b}\gamma,Q,PB\big)$. Further, $\dot{U}$ is non-positive, because $\|e\|\geq v_\mathrm{e}$ if $U\geq 0$, and $g_\mathrm{e}\geq\lambda_{\min}(Q)>0$. Observing $U\leq\lambda_{\max}(P)(\|e\|^2-{v_\mathrm{e}}^2)$ yields the desired statement. Let $t^\star\!=\!\text{min}\big\{t\!\geq\! 0:U\big|_t\!=\! 0\big\}$ if existent, and $t^\star=\infty$ otherwise. Given $\|e(0)\|\geq r_\mathrm{e}$, we have  $U\big|_{t=0}\geq 0$. That together with  $U\big|_t$ being continuous with respect to $t$ yields  $U\big|_t\geq 0$ for all  $t\in [0,t^\star)$. Applying the comparison principle to $\dot{U}\big|_t\leq-c_\mathrm{e} U\big|_t$ results in $(\ref{SlowestConvergence})$ for all $t\in [0,t^\star)$. Further, we have valid $\|e(t)\|\leq r_\mathrm{e}$ for all $t\in[t^\star\!,\infty)$, because the set $\Omega^\mathrm{e}=\big\{e\in\mathbb{R}^n:U\leq 0\big\}=\Omega_{\rho(0)}^\mathrm{e}$ is positively invariant and $\|e\|\leq r_\mathrm{e}$ holds for all $e\in\Omega^\mathrm{e}$.\vspace{-0.75mm}
\end{pf}

The properties of function $g\big(\cdot,Q,PB\big)$ listed in Lemma~\ref{lemma.EigenvalueLift} allow for some qualitative statements regarding the bound on the transient error~(\ref{SlowestConvergence}): Due to $\mathscr{b}>0$, there exists a nonempty interval $\mathcal{J}=(0,\sup\, J)$ such that the minimal eigenvalue $g(\alpha\mathscr{b}\gamma,Q,PB)$ in $r_\mathrm{e}$ and $c_\mathrm{e}$ is strictly increasing with $\alpha\in \mathcal{J}$ if and only if the inequality
\begin{align}
\label{EigenvalueHub}
    g\big(1,Q,PB\big)>\lambda_{\min}(Q)
\end{align}
holds and constant $\gamma$ is chosen in $(0,1)$. We readily obtain:
\begin{cor}\label{corollary:ConvergenceRate}\vspace{-0.25mm}    
Let $\mathscr{b}>0$ and evaluate $r_\mathrm{e}\big|_{\mu=0}$ in the bound on the transient error $(\ref{SlowestConvergence})$ for a fixed $\gamma\in(0,1)$. Then, the convergence rate $c_\mathrm{e}$ of $(\ref{SlowestConvergence})$ is strictly increasing with scaling $\alpha\in \mathcal{J}$ of the gain of the static update law $(\ref{StaticUpdateLaw})$, where $\mathcal{J}=(0,\sup\, \mathcal{J})$ is some nonempty interval, if and only if $(\ref{EigenvalueHub})$ holds. For $n\geq 2$, however, the increase in the convergence rate is subject to the upper bound $c_\mathrm{e}\leq \scriptstyle \frac{\lambda_{\max}(Q)}{\lambda_{\max}(P)}$, regardless of the size of $\alpha$. In contrast, the convergence rate stays at $ c_\mathrm{e }= \scriptstyle \frac{\lambda_{\min}(Q)}{\lambda_{\max}(P)}$ for all $\alpha>0$ if $(\ref{EigenvalueHub})$ does not apply.\vspace{-1mm}
\end{cor}
\begin{pf}
    Use Lemma~\ref{lemma.EigenvalueLift} except for the trivial case $n=1$.\vspace{-0.75mm}
\end{pf}

The inequality $(\ref{EigenvalueHub})$ can either be checked after $0<P\in\mathrm{S}^n$ for a given $0<Q\in\mathrm{S}^n$ has been calculated or, in order to ensure an increase in the convergence rate $c_\mathrm{e}$ with $\alpha$, it can serve as an additional condition for the solution of the Lyapunov equation~(\ref{NomStab}) besides $0<P,Q\in\mathrm{S}^n\!$. Section~\ref{part.GraphicalCriterion} shows a criterion for the feasibility of the latter approach.\linebreak The influence of the constant $\gamma\in[0,1)$ on the value $r_\mathrm{e}\big|_{\mu=0}$ has to be discussed in the light of the inequality $(\ref{EigenvalueHub})$. If $(\ref{EigenvalueHub})$ is not satisfied, i.e. if the convergence rate $c_\mathrm{e}$ is not affected by $\alpha$, then $r_\mathrm{e}\big|_{\mu=0}$ at the minimum point $\gamma=0$ is our\linebreak least conservative estimate of the residual error of $\|e\|$. In the case of $(\ref{EigenvalueHub})$ being fulfilled, the constant $\gamma$ represents a compromise between the attenuation of unstructured uncertainty and structured uncertainty in $r_\mathrm{e}\big|_{\mu=0}$, which we generally cannot optimally resolve without knowing $W$, $\eta^\star$ and $g(\cdot,Q,PB)$. However, specific knowledge of the function $g(\cdot,Q,PB)$ allows for choosing $\gamma$ in a suboptimal manner without requiring $W$ or $\eta^\star$, as shown below.

Suppose, $0<P,Q\in\mathrm{S}^n$ satisfy the inequality $(\ref{EigenvalueHub})$ and the function $g(\cdot,Q,PB)$ is differentiable on $(0,\mathscr{d})$ with some $\mathscr{d}>0$. If the gain of the static update law is scaled with
\begin{align}
    \label{LowerBoundAlpha}
    \alpha>\frac{\lambda_{\min}(Q)}{\mathscr{b}}\left(\lim_{\phi\searrow 0}\;{\frac{\partial}{\partial\phi}g\big(\phi,Q,PB\big)}\right)^{-1},
\end{align}
then the value $\Phi^\star=\mathrm{sup} \mathcal{M}$ with the set
\begin{align}
\label{phiStar}
\!\!\!\mathcal{M}\!=\!\left\{\!\Phi\!\in\!(0,\!\alpha\mathscr{b}]\!:\!\phi\!+\!\frac{g\big(\phi,Q,PB\big)}{\frac{\partial}{\partial\phi}g\big(\phi,Q,PB\big)}\!\leq\!\alpha\mathscr{b}\;\forall \phi\!\in\! (0,\!\Phi] \!\right\}\!\!\!\!
\end{align}
renders $\gamma^\star\!=\!\Phi^\star(\alpha\mathscr{b})^{-1}\!\in\!(0,1)$ suboptimal in the sense that\linebreak $r_\mathrm{e}\big|_{\substack{\gamma=\gamma^\star\\ \hspace{-1ex}\mu=0}}\!\leq\! r_\mathrm{e}\big|_{\substack{\gamma=0 \\ \mu=0}}$. Due to $\alpha\mathscr{b}\gamma^\star\!>\!0$, also $c_\mathrm{e}\big|_{\gamma=\gamma^*}\!>\!{\scriptstyle \frac{\lambda_{\min}(Q)}{\lambda_{\max}(P)}}$. \vspace{-2mm}
\begin{pf}
    Assuming $\Phi^\star=\alpha\mathscr{b}$ yields the contradiction $\lim_{\phi\nearrow\alpha\mathscr{b}}{ \scriptstyle \frac{g(\phi,Q,PB)}{\partial g(\phi,Q,PB)/\partial\phi}}=0$ to the fact that $g(\cdot,Q,PB)$ is Lipschitz continuous, i.e. $\partial g(\cdot,Q,PB)/ {\partial\phi}$ is e.b., since $g(\phi,Q,PB)\geq\lambda_{\min}(Q)>0$ for all $\phi\geq 0$. Further, $\mathcal{M}\!=\!\emptyset$ would be inconsistent with $(\ref{LowerBoundAlpha})$. Hence, $\Phi^\star\in(0,\alpha\mathscr{b})$ and $\gamma^\star\in(0,1)$ hold. We get $\forall \phi\in(0,\Phi^\star]: \phi {\scriptstyle \frac{\partial}{\partial\phi}}g(\phi,Q,PB)+g(\phi,Q,PB)\leq \alpha \mathscr{b} {\scriptstyle \frac{\partial}{\partial\phi}}g(\phi,Q,PB)$ from $(\ref{phiStar})$. Substituting $\phi=\alpha\mathscr{b}\gamma$ leads to $\forall \gamma\in(0,\gamma^\star]: {\scriptstyle \frac{\mathrm{d}}{\mathrm{d}\gamma}}(1-\gamma)g(\alpha\mathscr{b}\gamma,Q,PB)=\alpha\mathscr{b}(1-\gamma){\scriptstyle \frac{\partial}{\partial\phi}}g(\phi,Q,PB)\big|_{\phi=\alpha\mathscr{b}\gamma}-g(\alpha\mathscr{b}\gamma,Q,PB)\!\geq\! 0$, which in turn implies $g(0,Q,PB)\leq (1-\gamma^\star)g(\alpha\mathscr{b}\gamma^\star,Q,PB)$.\vspace{-0.75mm}
\end{pf}
In summary, we obtain: If the solution $0<P,Q\in\mathrm{S}^n$ of the Lyapunov equation~(\ref{NomStab}) fulfills the inequality $(\ref{EigenvalueHub})$ and the gain of the static update law is large enough to satisfy $(\ref{LowerBoundAlpha})$, then the norm of the transient error approaches our least conservative estimate $\inf_{\gamma\in[0,1)}\big\{r_\mathrm{e}\big|_{\mu=0}\big\}$ of its residual error with increased convergence rate $c_\mathrm{e} > \scriptstyle \frac{\lambda_{\min}(Q)}{\lambda_{\max}(P)}$.

\begin{rem}
    The derivations of $r_\mathrm{e}$ and $c_\mathrm{e}$ both neglect the influence of the term $-\tilde{W}^\top\!K^{-1} \tilde{W}\!=\!-\alpha^{-1}\tilde{W}^\top\! K_\mathrm{b}^{-1} \tilde{W}$ in $(\ref{EqdotVStaticUpdateLaw})$. Thus, for larger $\alpha$, i.e. higher gain of the static update law, $r_\mathrm{e}$ is a more accurate estimate of the actual ultimate bound, the actual residual error, and $c_\mathrm{e}$ is a closer approximation of the actual convergence rate of $\|e\|$.\vspace{-0.75mm}
\end{rem}

\vspace{-1.25mm}\subsection{Structured Uncertainty with Zero Residual Error for Adaptation with the Static Update Law}\label{part:ResidualError}\vspace{-1.25mm}
The residual error is an upper bound for the norm of the
error state in the steady state. With function $\beta$ in the structured uncertainty satisfying Assumption~\ref{assumption:WellPosedness}, zero residual error for adaptation with the static update law is not guaranteed---even in the absence of unstructured uncertainty. Here, we show that this drawback is removed if we further restrict the structured uncertainty. By referring to methods that quantify the robustness of multivariable systems, we calculate a linear-growth bound for the structured uncertainty such that the residual error is zero when no unstructured uncertainty is present in the closed loop. The studies by \cite{patel_quantitative_1980}, \cite{chen_improved_1994}, \cite{kim_comments_1995} let estimate the robustness against uncertainty in unspecified input channels, where \cite{kim_comments_1995} proposes the least conservative approach. We utilize the machinery in \cite{kim_comments_1995} in order to find a linear-growth bound for the structured uncertainty in the input channel of the error system $(\ref{ErrorDyn})$. We use the notations $(\ref{GainScaling})$ and $(\ref{LowerBoundBeta})$. Reorganizing $(\ref{EqdotVStaticUpdateLaw})$ with $\eta\equiv 0$ yields
\begin{align}
    \dot{V}\big|_{(\ref{StaticUpdateLaw})}\!=\!&-\!e^\top \big(Q\!+\!2(\beta^\top K\beta) PBB^\top P\big)e\!+\!2e^\top PB W^\top \beta\nonumber\\
    \leq\!& -\!e^\top \big(Q\!+\!2\alpha\mathscr{b}PBB^\top P\big)e\!+\!2e^\top P B W^\top \beta\nonumber\\
    \leq\!& -\!e^\top{\big(Q\!+\!(2\alpha\mathscr{b}- \varepsilon)PBB^\top P\big)}e\!+\!\frac{1}{ \varepsilon}\big(W^\top \beta\big)^2\label{Ineq1}
\end{align} 
for any $\varepsilon>0$. Subsequently, we allow for negative values in the first argument of function $g$. Consider the set
\begin{align}
\mathcal{E}=&\big\{ \varepsilon>0:g\big(2\alpha\mathscr{b}- \varepsilon,Q,PB\big)>0\big\}.
\end{align}
By restricting the constant $ \varepsilon>0$ to $ \varepsilon\in\mathcal{E}$, we get 
\begin{align}
\label{Ineq2}
    \dot{V}\big|_{(\ref{StaticUpdateLaw})}\leq&-g\big(2\alpha\mathscr{b}- \varepsilon,Q,PB\big)\|e\|^2+\frac{1}{ \varepsilon}\big|W^\top \beta\big|^2
\end{align}
from $(\ref{Ineq1})$. This leads to the next result.
\begin{thm}\label{theorem:ZeroResidualError} Consider the error dynamics $(\ref{ErrorDyn})$. Let Assumption $\ref{assumption:WellPosedness}$ apply and suppose, $0<P,Q\in\mathrm{S}^n$ satisfy the Lyapunov equation $(\ref{NomStab})$. Then, the adaptive control $(\ref{ControlMRAC})$ and the static update law $(\ref{StaticUpdateLaw})$ with gain $0<K\in\mathrm{S}^n$ render the origin of $e$ GAS, i.e. $\lim_{t\to\infty}\|e(t)\|=0$, if $\eta\equiv 0$ and
\begin{align}
\label{Robustness}
    \!\!\!\forall x\!\in\!\mathbb{R}^n\!\setminus\!\{0\}\!:\!\frac{\big|W^\top\!\beta(x)\big|}{\|x\|}\!<\!\underbrace{\max_{ \varepsilon\in\mathcal{E}}\!\sqrt{ \varepsilon g\big(2\alpha\mathscr{b}\!-\! \varepsilon,Q,PB\big)}}_{=\tau}.\!
\end{align}
\end{thm}

\begin{pf}
    First, we show the existence of the maximum in $(\ref{Robustness})$. The matrices $P,Q\in\mathrm{S}^n$ are positive definite. Thus, $\mathcal{E}\neq\emptyset$ since $2\alpha\mathscr{b}+\lambda_{\min}(Q)/\big(2\|PB\|^2\big)\in\mathcal{E}$. We have $\varepsilon g\big(2\alpha\mathscr{b}- \varepsilon,Q,PB\big)=0$ for $ \varepsilon\in\partial \mathcal{E}=\{0,\sup \mathcal{E}\}$ since $g\big(2\alpha\mathscr{b}-\sup\mathcal{E},Q,PB\big)=0$ by design, $\sup \mathcal{E}<\infty$ and $g\big(2\alpha\mathscr{b},Q,PB\big)\leq \lambda_{\max}(Q)<\infty$. Now, the maximum in question exists, because $ \varepsilon g\big(2\alpha\mathscr{b}- \varepsilon,Q,PB\big)$ is continuous with respect to $ \varepsilon\in\mathcal{E}\cup\partial\mathcal{E}$, and $ \varepsilon g\big(2\alpha\mathscr{b}- \varepsilon,Q,PB\big)>0$ holds for $ \varepsilon=2\alpha\mathscr{b}+\lambda_{\min}(Q)/\big(2\|PB\|^2\big)\in \mathcal{E}\setminus \partial\mathcal{E}$. Next, with $(\ref{Robustness})$, inequality $(\ref{Ineq2})$ yields $\dot{V}\big|_{(\ref{StaticUpdateLaw})}<0$ for all $e\neq 0$. GAS of $e=0$ results from usual Lyapunov theory.\vspace{-0.75mm}
\end{pf}

If the function $\beta$ in addition to satisfying Assumption~\ref{assumption:WellPosedness} is bounded from below with $\mathscr{b}>0$, then we can evaluate the linear-growth bound $\tau$ in a qualitative manner: Due to $0<2\alpha\mathscr{b}\in\mathcal{E}$, we have $\tau\geq \big(2\alpha\mathscr{b}\lambda_{\min}(Q)\big)\!^{1/2}$, i.e. then the linear-growth bound $\tau$ can be arbitrarily increased by selecting the gain of the static update law sufficiently large. 

\begin{rem}
Function $\beta$ in the adaptive control $(\ref{ControlMRAC})$ and in the static update law $(\ref{StaticUpdateLaw})$ can always be replaced by the augmented version  $(\beta^\top\ 1)^\top$, which gives $\mathscr{b}\geq \lambda_{\min}(K_\mathrm{b})>0$ and thereby renders the analysis of the ultimate bound and the transient behavior in Section~\ref{part:StaticUpdateLaw} applicable. The structured uncertainty that acts on the error dynamics $(\ref{ErrorDyn})$ remains unchanged by the corresponding augmentation $(W^\top\ 0)^\top$ of the weights $W$. Moreover, the linear-growth bound $\tau$ with respect to the augmented structured uncertainty can be arbitrarily increased via selecting sufficiently large $\alpha>0$, while augmenting the weights is leaving the term $\big|W^\top\beta(x)\big|/\|x\|$ in $(\ref{Robustness})$ unchanged for all $x\in\mathbb{R}^n\setminus\{0\}$.
\end{rem}

\vspace{-1.25mm}\subsection{Criterion based on the KYP Lemma}\label{part.GraphicalCriterion}\vspace{-1.25mm}
For adaptation with the static update law, the validity of inequality~(\ref{EigenvalueHub}) determines whether or not the bound on the transient error~(\ref{SlowestConvergence}) decays faster when increasing the gain of the update law, see Corollary~\ref{corollary:ConvergenceRate}. Thus, it is desirable to find a criterion for the existence of $0<P,Q\in\mathrm{S}^n$ satisfying $(\ref{EigenvalueHub})$ and the Lyapunov equation $(\ref{NomStab})$. We approach this problem by aiming at enforcing  $0<P,Q\in\mathrm{S}^n$, the Lyapunov equation~(\ref{NomStab}) and  $PB=v$ for a suitably selected\linebreak $v\in\mathbb{R}^n$. This simplifies the inequality~(\ref{EigenvalueHub}) to $g(1,Q,v)>\lambda_{\min}(Q)$, enabling the proof of the criterion in Theorem~\ref{theorem:GraphicalCriterion}.

\begin{thm}\label{theorem:GraphicalCriterion}
    Let $A\in\mathbb{R}^{n\times n}$ Hurwitz, $B,v\in\mathbb{R}^n$, $(A,B)$ controllable, scalar $\kappa>0$, $f:\mathbb{R}\times[0,\infty)\times[0,\infty)\times\mathbb{R}^n\to\mathbb{R}$ and $G:\mathrm{j}\mathbb{R}\to\mathbb{C}^n$ be functions with
    \begin{align}
        \!\!\!G(\mathrm{j}\omega)\!=&\big(\mathrm{j}\omega\mathrm{I}_n\!-\!A\big)\!^{\,-1}B,\\
         \!\!\!f\!\big(\omega,\!\kappa,\!\vartheta,\!v\big)\!=&2\,\mathrm{Re}\big\{v\!^\top\!G(\mathrm{j}\omega)\big\}\!+\!\kappa\big|v\!^\top\!G(\mathrm{j}\omega)\big|^2\!\!-\!\vartheta\big\|G(\mathrm{j}\omega)\big\|^2\!\!.\!\!
    \end{align}
    Further, consider the linear matrix inequality (LMI)
    \begin{align}
    \label{LMI}
        A^\top P+PA+\big(\psi+\kappa\|v\|^2\big)\mathrm{I}_n-\kappa vv^\top\leq&0,\quad\!\! PB=v
    \end{align}
    in dependence of the scalar $\psi>0$. The set
    \begin{align}
    \Psi=&\big\{\psi>0:\text{The LMI }(\ref{LMI})\text{ for }P\in\mathrm{S}^n\text{ is feasible} \big\}
    \end{align}
    is nonempty if and only if
    \begin{align}
    \label{Criterion}
      &\forall \omega\in\mathbb{R}:f\big(\omega,\kappa,\kappa\|v\|^2,v\big)>0,\\
         \label{Limes}
        &\lim_{\omega\to\infty}\big\|G(\mathrm{j}\omega)\big\|^{-2}f\big(\omega,\kappa,\kappa\|v\|^2,v\big)>0.
    \end{align}
    If $\Psi$ is nonempty, then any solution $P\!\in\!\mathrm{S}^n$ of the LMI $(\ref{LMI})$ with $\psi\in\Psi$ satisfies $P>0$ and $Q=-\big(A^\top P+PA\big)>0$. Moreover, if the criterion
    \begin{align}
    \label{ConditionEigenvalueHub}
    \exists\omega\in\mathbb{R}:f\big(\omega,0,\psi+\kappa\|v\|^2,v\big)<0
    \end{align}
    is fulfilled, then the inequality $(\ref{EigenvalueHub})$ is satisfied.\vspace{-0.75mm}
\end{thm}

\begin{pf}
    By means of Lemma~\ref{lemma:BlockMatrices} from the Appendix, we have $\psi\in\Psi$ if and only if 
    \begin{align}
    \label{MatrixIneqaulityKYP}
        \!\!\!\!\begin{pmatrix}
            A^\top P+PA+\big(\psi+\kappa \|v\|^2\big)\mathrm{I}_{n}-\kappa vv^\top&PB-v\\\big(PB-v\big)^\top&0
        \end{pmatrix}\leq& 0.\!
    \end{align}
    The KYP lemma, as it is stated in \cite{rantzer_kalmanyakubovichpopov_1996}, together with $A$ Hurwitz and the controllability of $(A,B)$ gives the following: There exists $P\!\in\!\mathrm{S}^n$ satisfying~(\ref{MatrixIneqaulityKYP}) if and only if\vspace{-0.75mm}
    \begin{align}
    \label{FrequencyIneqaulityKYP}
        \!\!\!\!\begin{pmatrix}
            G(\mathrm{j}\omega)\\1
        \end{pmatrix}^{\!\mathrm{h}}\begin{pmatrix}
            \big(\psi+\kappa \|v\|^2\big)\mathrm{I}_{n}-\kappa vv^\top&-v\\-v^\top&0
        \end{pmatrix}\begin{pmatrix}
            G(\mathrm{j}\omega)\\1
        \end{pmatrix}\leq& 0\!
    \end{align}
    for all $\omega\in\mathbb{R}\cup\{\infty\}$, which in turn is equivalent to $\forall \omega\in\mathbb{R}:f\big(\omega,\kappa,\psi+\kappa\|v\|^2,v\big)\geq 0$. Thus, the set $\Psi$ can be written as $\Psi=\{\psi>0:f\big(\omega,\kappa,\psi+\kappa\|v\|^2,v\big)\geq\! 0\;\forall \omega\in\mathbb{R}\}$. Observing $f\big(\omega,\kappa,\psi+\kappa\|v\|^2,\!v\big)= f\big(\omega,\kappa,\kappa\|v\|^2,\!v\big)-\psi\|G(\mathrm{j}\omega)\|^2$
    for all $\omega\in\mathbb{R}$ yields the equivalence $\psi\in\Psi$ if and only if $0<\psi\leq \|G(\mathrm{j}\omega)\|^{-2}f\big(\omega,\kappa,\kappa\|v\|^2,v\big)$ for all $\omega\in\mathbb{R}$, holding, if and only if $0<\psi\leq \mathrm{inf}_{\omega\in\mathbb{R}}\|G(\mathrm{j}\omega)\|^{-2}f\big(\omega,\kappa,\kappa\|v\|^2,v\big)$. The transfer function $G$ has no poles on the imaginary axis, since $A$ is Hurwitz. Therefore, $\mathrm{inf}_{\omega\in\mathbb{R}}\|G(\mathrm{j}\omega)\|^{-2}f\big(\omega,\kappa,\kappa\|v\|^2,v\big)>0$, i.e. $\Psi\neq\emptyset$, if and only if $(\ref{Criterion})$ and $(\ref{Limes})$ hold. For the remainder of this proof, consider $(\ref{Criterion})$ and $(\ref{Limes})$ to be satisfied and take $\psi\in\Psi$. The solution $P\in\mathrm{S}^n$ of the LMI $(\ref{LMI})$ admits $-Q+\big(\psi+\kappa\|v\|^2\big)\mathrm{I}_n-\kappa vv^\top\leq 0$. We have $\kappa\|v\|^2\mathrm{I}_n-\kappa vv^\top\geq 0$ and therefore $Q>0$ holds since $\psi>0$. Additionally, $P\geq 0$ holds by virtue of Lemma~$1.10.1$ in \cite{leonov_frequency-domain_1996} since $A$ is Hurwitz. In order to verify $P>0$ by means of a contradiction, suppose $\mathrm{det}(P)=0$. Then, there exists $0\neq w\in\mathbb{R}^n$ with $Pw=0$. Due to $P\in\mathrm{S}^n$, we receive the contradiction $0=2w^\top A^\top Pw=-w^\top Qw$ to $w^\top Qw>0$. Thus, the matrix $P$ has no eigenvalues at zero, hence $P>0$ because all eigenvalues of $P$ are non-negative from the outset. Furthermore, the solution admits $\lambda_{\min}\big(Q+\kappa vv^\top\big)\geq\psi+\kappa\|v\|^2$. Assume, $(\ref{ConditionEigenvalueHub})$ holds. The KYP lemma in \cite{rantzer_kalmanyakubovichpopov_1996} as well as Lemma~\ref{lemma:BlockMatrices} both state equivalences. Using these equivalences on $(\ref{ConditionEigenvalueHub})$ while negating their statements results in  $A^\top\bar{P}+\bar{P}A+\big(\psi+\kappa\|v\|^2\big)\mathrm{I}_n\nleq 0$, i.e. $\lambda_{\min}\big(-(A^\top\bar{P}+\bar{P}A)\big)<\psi+\kappa\|v\|^2$, or $\bar{P}B\neq v$ for all $\bar{P}\in\mathrm{S}^n$. The choice $\bar{P}=P$ shows $\lambda_{\min} (Q)<\psi+\kappa vv^\top$ since $PB=v$. The desired inequality $g(1,Q,v)>\lambda_{\min}(Q)$ follows from $\lambda_{\min}\big(Q+\kappa vv^\top\big)\geq \psi+\kappa\|v\|^2\!>\!\lambda_{\min}(Q)$ via negation of the statements in the equivalence~(\ref{EigenvalueLiftProof}) since $\kappa>0$ and $g(1,Q,v)\geq \lambda_{\min}(Q)$.\vspace{-0.75mm}
\end{pf}

There are two obstacles in Theorem~\ref{theorem:GraphicalCriterion} that hinder its application for finding a solution $0<P,Q\in\mathrm{S}^n$ of the Lyapunov equation $(\ref{DefLyapunovFunction})$ that satisfies the inequality $(\ref{EigenvalueHub})$: First, one has to find $v\in\mathbb{R}^n$ such that $(\ref{Criterion})$ and $(\ref{Limes})$ hold. Second, since the set $\Psi$ is not explicitly known, one has to carefully tune the slack variable $\psi>0$ so that the LMI $(\ref{LMI})$ is feasible and the criterion $(\ref{ConditionEigenvalueHub})$ is satisfied. These difficulties are mitigated in Corollary~\ref{corollary:GraphicalCriterion} by means of a \mbox{$\psi$-free} criterion for the validity of the inequality $(\ref{EigenvalueHub})$ for any solution $P\in\mathrm{S}^n$ of the LMI 
\begin{align}
    \label{strictLMI}
        A^\top P+PA+\kappa\|v\|^2\mathrm{I}_n-\kappa vv^\top<&0,\quad PB=v.
\end{align}

\begin{defn} \citep{ioannou_frequency_1987} Let $A\in\mathbb{R}^{n\times n}$ and $B,v\in\mathbb{R}^n$. Transfer function $H(s)=v^\top\big(s\mathrm{I}_n-A\big)\!^{\,-1}B$ with $s\in\mathbb{C}$ is called strictly positive real (SPR) if $H(s)$ has no poles in $\mathrm{Re}\{s\}\geq 0$, $\forall \omega\in\mathbb{R}: \mathrm{Re}\big\{H(\mathrm{j}\omega)\big\}>0$ and $\mathrm{lim}_{\omega\to\infty}\omega^2\,\mathrm{Re}\{H(\mathrm{j}\omega)\}>0$.
\end{defn}

\begin{cor}\label{corollary:GraphicalCriterion}
Let $A\in\mathbb{R}^{n\times n}$ Hurwitz, $B,v\in\mathbb{R}^n$ and $(A,B)$ controllable. If the transfer function $H(s)=v^\top\big(s\mathrm{I}_n-A\big)\!^{\,-1}B$ is SPR, then there exists a nonempty interval $\mathcal{K}=(0,\sup\mathcal{K})$ such that, for each $\kappa\in\mathcal{K}$, the LMI $(\ref{strictLMI})$ for $P\in\mathrm{S}^n$ is feasible. For any scalar $\varrho\in[0,1)$, an admissible upper bound of $\mathcal{K}$ is given by $\sup\mathcal{K}=\mathrm{inf}_{\omega\in\mathbb{R}}\big\{2\,\mathrm{Re}\{H(\mathrm{j}\omega)\}/\big(\|v\|^2\|G(\mathrm{j}\omega)\|^2-\varrho |H(\mathrm{j}\omega)|^2\big)\big\}$. Provided that the LMI $(\ref{strictLMI})$ for $\kappa>0$ has a solution $P\in\mathrm{S}^n$, then $P>0$ and $Q=-\big(A^\top P+PA\big)>0$. Moreover, the inequality $(\ref{EigenvalueHub})$ is satisfied if the \mbox{$\psi$-free} criterion 
 \begin{align}
     \label{psiFreeConditionEigenvalueHub}
    \exists\omega\in\mathbb{R}:f\big(\omega,0,\kappa\|v\|^2,v\big)\leq0
 \end{align}
 is fulfilled.\vspace{-0.75mm}
\end{cor}

\begin{pf}
    Suppose, $H(s)$ is SPR. First, we show that there is a nonempty interval $\mathcal{K}=(0,\sup\mathcal{K})$ such that $(\ref{Criterion})$ and $(\ref{Limes})$ hold for all $\kappa\in\mathcal{K}$. The inequality $(\ref{Criterion})$ is satisfied if $\kappa<\bar{\kappa}(\omega)$ for all $\omega\in\mathbb{R}$ with $\bar{\kappa}(\omega)=2\,\mathrm{Re}\{H(\mathrm{j}\omega)\}/\big(\|v\|^2\|G(\mathrm{j}\omega)\|^2-\varrho|H(\mathrm{j}\omega)|^2\big)$, holding, if and only if $\kappa<\mathrm{inf}_{\omega\in\mathbb{R}}\bar{\kappa}(\omega)$. We have $\mathrm{Re}\{H(\mathrm{j}\omega)\}>0$ and there exists $b^\star\in(0,\infty)$ such that $\|v\|^2\|G(\mathrm{j}\omega)\|^2-\varrho|H(\mathrm{j}\omega)|^2=G(\mathrm{j}\omega)^\mathrm{h}\big(\|v\|^2\mathrm{I}_n-\varrho vv^\top\big)G(\mathrm{j}\omega)\in(0,b^\star]$ for all $\omega\in\mathbb{R}$ since $H(s)$ is SPR, $\|v\|^2\mathrm{I}_n-\varrho vv^\top>0$ and $A$ is Hurwitz, i.e. $0\!<\!\|G(\mathrm{j}\omega)\|\!<\!\infty$ for all $\omega\!\in\!\mathbb{R}$. Hence, $\mathrm{inf}_{\omega\in\mathbb{R}}\bar{\kappa}(\omega)\!>\!0$ if and only if $\mathrm{lim}_{\omega\to\infty}\bar{\kappa}(\omega)>0$, where the latter is true since $\|v\|^2\|G(\mathrm{j}\omega)\|^2-\varrho|H(\mathrm{j}\omega)|^2$ is a rational function in $\omega$ with a relative degree of at least two and, due to $H(s)$ being SPR, we have $\mathrm{lim}_{\omega\to\infty}\omega^2\, \mathrm{Re}\{H(\mathrm{j}\omega)\}>0$. In short, the inequality $(\ref{Criterion})$ holds for all $\kappa\in(0,\mathrm{inf}_{\omega\in\mathbb{R}}\bar{\kappa}(\omega))\neq\emptyset$. Next, the inequality $(\ref{Limes})$ can be rewritten as $\mathrm{lim}_{\omega\to\infty}\big(\bar{\kappa}(\omega)-\kappa\big)\big(\|v\|^2-\varrho\|G(\mathrm{j}\omega)\|^{-2}|H(\mathrm{j}\omega)|^2\big) + \kappa(1-\varrho)\|G(\mathrm{j}\omega)\|^{-2}|H(\mathrm{j}\omega)|^2>0$, which is satisfied for all $\kappa\in(0,\mathrm{inf}_{\omega\in\mathbb{R}}\bar{\kappa}(\omega))$  because $\varrho<1$ and $\|G(\mathrm{j}\omega)\|^{-2}|H(\mathrm{j}\omega)|^2=\|G(\mathrm{j}\omega)\|^{-2} G(\mathrm{j}\omega)^\mathrm{h} vv^\top G(\mathrm{j}\omega)\leq\|v\|^2$ for all $\omega\in\mathbb{R}$ as a result of $\lambda_{\mathrm{max}}(vv^\top)=\|v\|^2$. Now, by virtue of Theorem~\ref{theorem:GraphicalCriterion}, $\Psi\neq\emptyset$ holds for each $\kappa\in\mathcal{K}$. Whenever $\psi\in\Psi\neq\emptyset$, then the LMI $(\ref{strictLMI})$ is feasible for the solution $P\in\mathrm{S}^n$ of the LMI $(\ref{LMI})$. Further, by construction of the set $\Psi$, for every solution $P\in\mathrm{S}^n$ of the LMI $(\ref{strictLMI})$ with $\kappa>0$, there exist $\psi\in\Psi$ such that $P$ is also a solution of the LMI $(\ref{LMI})$, e.g. $\psi=\lambda_{\min}\big(Q+\kappa vv^\top-\kappa\|v\|^2 \mathrm{I}_n\big)>0$. Thereby, the remaining statements of Corollary~\ref{corollary:GraphicalCriterion} follow from Theorem~\ref{theorem:GraphicalCriterion}, because if the \mbox{$\psi$-free} criterion $(\ref{psiFreeConditionEigenvalueHub})$ is satisfied, then the criterion $(\ref{ConditionEigenvalueHub})$  for all $\psi\in\Psi$ is also satisfied since $\|G(\mathrm{j}\omega)\|>0$ for all $\omega\in\mathbb{R}$.\vspace{-0.75mm}
\end{pf}

The applicability of Corollary~\ref{corollary:GraphicalCriterion} to ensure that inequality $(\ref{EigenvalueHub})$ holds also depends  on properly tuning a slack variable, namely $\varrho$. However, $\varrho\in[0,1)$ only needs to be tuned so that the $\psi$-free criterion $(\ref{psiFreeConditionEigenvalueHub})$ for a $\kappa\in\mathcal{K}$ is feasible. That is, unlike the slack variable  $\psi$ in Theorem~\ref{theorem:GraphicalCriterion}, the tuning of $\varrho$ does not depend on the feasibility of an LMI.



\vspace{-1.25mm}\section{Illustrative Example: Second Order Plant without Transfer Zeros} \label{section:SIMULATION}\vspace{-1.25mm}
We consider generic second order dynamics with natural frequency $\omega_0=1$ and damping $\zeta=2^{-1/2}$ for the error system $(\ref{ErrorDyn})$ to represent a mass moving along a trajectory. The system matrices can be found in the Appendix. This error system is able to oscillate but the frequency response of the position error exhibits no resonance peak. We compare the closed-loop performance for the update laws\vspace{-.75mm}
\begin{align}
    \label{PIUpdateLaw_Sim}
    \text{PI: }&\hat{W}\!=\!K_\mathrm{PI}\beta B^\top P e\!+\!\Gamma_\mathrm{PI}\int_0^{t}\beta B^\top P e\!-\!\Sigma_\mathrm{PI}\hat{W}\mathrm{d}t,\\
    \label{StaticUpdateLaw_Sim}
    \text{static: }&\hat{W}\!=\!K_\mathrm{P}\beta B^\top P e
\end{align}
with $0\!<\!\Gamma_\mathrm{PI},\Sigma_\mathrm{PI},K_\mathrm{P}\!\in\! \mathrm{S}^{n_\beta}$, $K_\mathrm{PI}\!\in\! \mathrm{S}^{n_\beta}$ and $0\!\leq\! K_\mathrm{PI}\!<\!4\Sigma_\mathrm{PI}^{-1}$. Applying the $\psi$-free criterion $(\ref{psiFreeConditionEigenvalueHub})$ in Fig.~\ref{figure:GraphicalCriterion} leads to \mbox{$0\!<\!P\!\in\!\mathrm{S}^2$} in the Appendix that ensures an increase in the convergence rate $c_\mathrm{e}$ through linear scaling of the gain $K_\mathrm{P}$. The adaptation rate and the $\sigma$-modification are set to $\Gamma_\mathrm{PI}\!=\!2\,\mathrm{I}_{n_\beta}$ and $\Sigma_\mathrm{PI}\!=\!0.2\,\mathrm{I}_{n_\beta}$, respectively.

\begin{figure}
\centerline{\includegraphics[scale=0.5925]{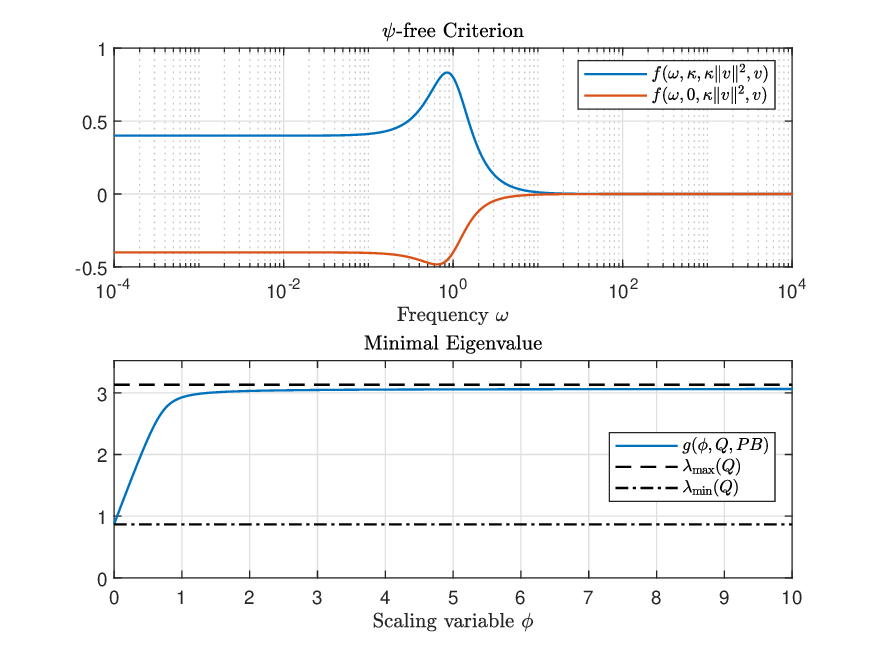}}
\vspace{-3mm}
\caption{We use $v\!=\!(1\ \ 2/2^{1/2})^\top$ such that $H(s)$ is SPR.\linebreak  The slack variables are set to $\varrho\!=\!0.75$ and $\kappa\!=\!0.9\,\sup\mathcal{K}$. The increase in the minimal eigenvalue $g(\phi,Q,PB)$ with $\phi\!\geq\! 0$ shown in the lower figure is a consequence of Corollary~\ref{corollary:GraphicalCriterion} since the \mbox{$\psi$-free} criterion~(\ref{psiFreeConditionEigenvalueHub}) is fulfilled, compare the red graph. The solution $0\!<\!P\!\in\!\mathrm{S}^2$ of the associated LMI~(\ref{strictLMI}) can be found in the Appendix.}
\label{figure:GraphicalCriterion}\vspace{-.25mm}
\end{figure}

\vspace{-1.25mm}\subsection{Simulation}
\begin{figure}
\centerline{\includegraphics[scale=0.5925]{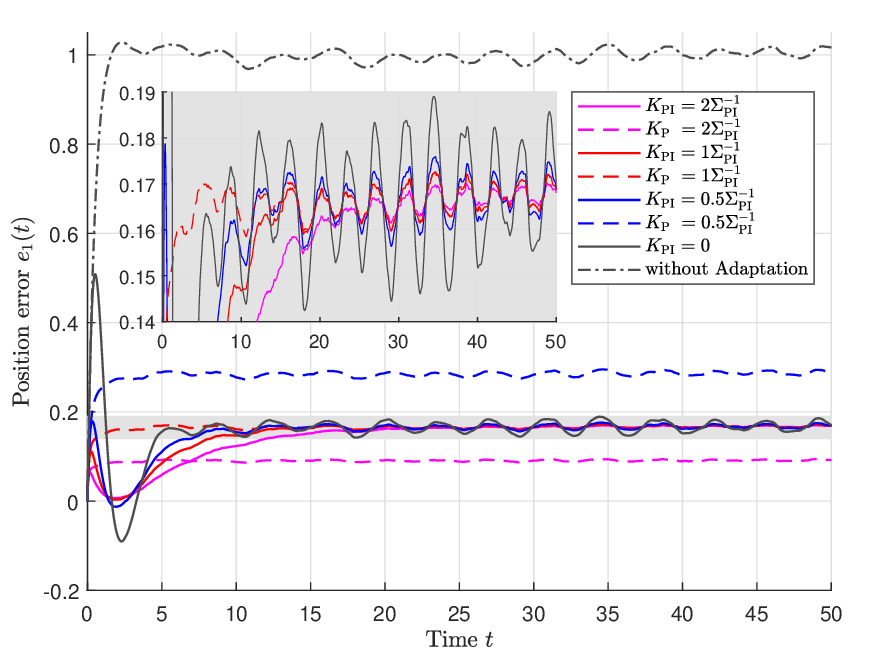}}
\vspace{-3mm}
\caption{Position error $e_1$ after an initial error in the velocity with different gains of the static update law (dashed) and the feedthrough in the PI update law (solid).}
\label{figure:PositionError}\vspace{-1mm}
\end{figure}

\vspace{-1.25mm}
The structured uncertainty is set to $\beta(e)^\top\!=\!(1\ -e_2)^\top$ where $W^\top\!=\!(1\ 1)^\top$ represents a constant offset and viscous friction. The unstructured uncertainty is set to $\eta(t)\!=\!\eta_\mathrm{w}(t)+0.05\,\mathrm{sin}(\omega_\mathrm{r}t)$ with $\omega_\mathrm{r}\!=\!1.7179$ and some band-limited white noise $\eta_\mathrm{w}(t)\!\in\!\mathbb{R}$. The simulation results in Fig. $\ref{figure:PositionError}$ show: First, adding feedthrough to the PI update law, or using the static update law, both yield significantly less oscillations of the position error. Second, for $K_\mathrm{P}\!=\!K_\mathrm{PI}\!=\!\Sigma_\mathrm{PI}^{-1}$, the steady-state behavior is indistinguishable between the static update law and the PI update law as, in this case, one update law tends to the other for large times. Third, increasing the steady-state noise rejection by larger gain of the feedthrough in the PI update law comes at the cost of slowly decaying transients. Adaptation with the static update law yields faster decaying transients than with the PI update law. Last, removing the integration from the PI update law by using the static update law with $K_\mathrm{P}\!=\!K_\mathrm{PI}$ decreases the remaining position error if and only if $K_\mathrm{PI}\!>\!\Sigma_\mathrm{PI}^{-1}$.

\vspace{-1.25mm}\subsection{Adaptation acting as a Disturbance Observer}
\begin{figure}
\centerline{\includegraphics[scale=0.5925]{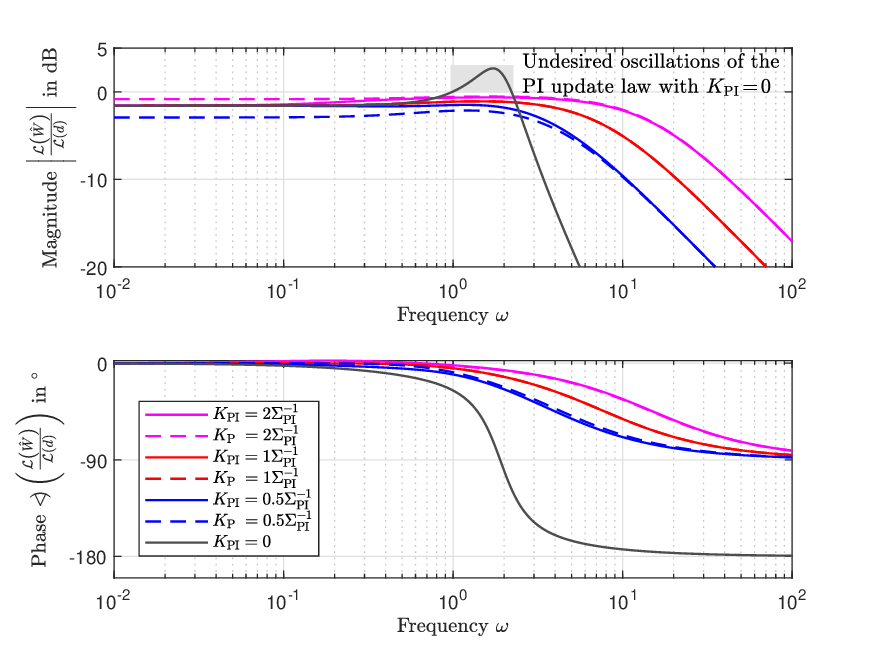}}
\vspace{-3mm}
\caption{Sensitivity of $\hat{W}(t)$ against the input disturbance $d(t)=W+\eta(t)$ for the purely linear case $\beta\equiv 1$ with different gains of the static update law (dashed) and the feedthrough in the PI update law (solid).}
\label{figure:DisturbanceObserver}
\end{figure}

\vspace{-1.25mm}
We consider the simplifying case $\beta\!\equiv\! 1$, for which all elements in the closed loop are linear such that we can evaluate the frequency-domain behavior of the adaptation.~Let
\begin{align}
    &\text{PI: }R\!=\!{\frac{K_\mathrm{PI}+s^{-1}\Gamma_\mathrm{PI}}{1+s^{-1}\Gamma_\mathrm{PI}\Sigma_\mathrm{PI}}}B^\top P,\; \text{static: }R\!=\!K_\mathrm{P}B^\top P.
\end{align}
Ignoring initial conditions, we get $\mathcal{L}(\hat{W})\!=\!R\mathcal{L}(e)$ and~thus
\begin{align}
    {\textstyle\frac{\mathcal{L}(\hat{W})}{\mathcal{L}(W+\eta)}}{=\!R\big(s\mathrm{I}_n-A+BR\big)\!^{\,-1}B}
\end{align}
from~(\ref{ErrorDyn}),~(\ref{ControlMRAC}). The calculated weight ${\hat{W}}(t)\!\in\!\mathbb{R}$ is acting as a disturbance observer of the remaining input noise $d(t)\!=\!W\!+\!\eta(t)$ in the following sense: If the sensitivity ${\scriptstyle\frac{\mathcal{L}(\hat{W})}{\mathcal{L}(d)}}|_{s=\mathrm{j}\omega}$ at a given frequency $\omega\!\in\!\mathbb{R}$ would be equal to one (i.e. $0\mathrm{dB}$, $0^\circ$), then sinusoidal input noise with this frequency would have no influence on the steady state of the error state, as it would be fully rejected by the adaptive control~(\ref{ControlMRAC}). We can deduce the following from Fig. \ref{figure:DisturbanceObserver}:~Undesired high-frequency oscillations of the disturbance observer for $K_\mathrm{PI}\!=\!0$, which result from high adaptation rates $\Gamma_\mathrm{PI}$ paired with insufficient damping $\Sigma_\mathrm{PI}$ of the PI update law, are suppressed by means of $K_\mathrm{PI}\!>\!0$ or the use of the static update law. Both,\linebreak adding feedthrough to the PI update law or using the static update law, broadens the bandwidth of the disturbance observer into higher frequencies, which comes at the cost of increased high-frequency control effort but causes better suppression of high-frequency input noise compared to the PI update law without feedthrough. For \mbox{$K_\mathrm{P}\!=\!K_\mathrm{PI}$}, the disturbance observer is only dependent on the choice of update law for excitation in the low-frequency band. Removing the integration from the PI update law, i.e. using the static update law with $K_\mathrm{P}\!=\!K_\mathrm{PI}$, yields increased precision of the disturbance observer if and only if \mbox{$K_\mathrm{PI}\!>\!\Sigma_\mathrm{PI}^{-1}$}, i.e. if the \mbox{$\sigma$-modification} is injecting too much damping into the PI update law. This underpins the results of the simulation. Removing the integration from an overly damped PI update law requires more low-frequency control effort. Note that the low-frequency behavior of the PI update law is independent of $K_\mathrm{PI}$. Further, both update laws for $K_\mathrm{P}\!=\!K_\mathrm{PI}\!=\!\Sigma_\mathrm{PI}^{-1}$ exhibit the same behavior in the steady state, i.e. the low-frequency behavior of the PI update law is also independent of the adaptation rate $\Gamma_\mathrm{PI}$. Regardless the magnitude of $K_\mathrm{PI}$, there are therefore only two options for increasing the suppression of low-frequency input noise, which is achieved by adaptation with the PI update law:\linebreak Decreasing $\Sigma_\mathrm{PI}$, or using the static update law with \linebreak$K_\mathrm{P}\!>\!\Sigma_\mathrm{PI}^{-1}$ instead, where the bandwidth of the disturbance observer then depends on the magnitude of $K_\mathrm{P}.$

\vspace{-1.4mm}
\vspace{-0.5mm}\section{Conclusion} \label{section:CONCLUSION}\vspace{-2mm}
We have simplified the indirect adaptive control approach by removing the integration from the update law. The analysis of the closed loop shows that the error state is ultimately bounded and that the ultimate bound can be forced arbitrarily close to the origin by increasing the gain of the static update law. Further, if we impose an additional condition on the solution of the Lyapunov equation besides nominal stability, we can guarantee that our bound on the transient error decays faster when increasing the gain of the static update law. This condition needs to be considered because of the semi-definite structure in the derivative of the Lyapunov function. We have presented a criterion for its fulfillment. In addition to the observation that adaptation with the static update law reduces undesired oscillations, for structured uncertainty satisfying a linear-growth bound, it forces the error state to its origin. In the purely linear case for a second order plant, the static update law compared to the PI update law provides better suppression of low-frequency input noise if and only if $K_\mathrm{P}\!>\!\Sigma_\mathrm{PI}^{-1}\!$, and it broadens the bandwidth of the adaptation into higher frequencies if and only if $K_\mathrm{P}\!>\!K_\mathrm{PI}$. 


\appendix
\vspace{-1.25mm}\section*{Appendix}\vspace{-1.25mm}

\begin{lem}\label{lemma:MinimalEigenvalue}
Let $Q\in\mathrm{S}^n$ and $w\in\mathbb{R}^{n\times m}$. The equality $g(1,Q,w)=\lambda_{\min}\left(Q\right)$ is valid if and only if there is  $r\in\mathbb{R}^n$, $r\neq 0$ with $Qr=\lambda_{\min}\left(Q\right) r$ such that $w^\top r=0$.\vspace{-0.75mm}
\end{lem}
\begin{pf}
    Suppose, there exists $p\in\mathbb{R}^n$, $\|p\|=1$ with $p^\top\left(Q+ww^\top\right)p=\lambda_{\min}\left(Q\right)$. Rearranging yields $\lambda_{\min}(Q)-p^\top Q p=\|w^\top p\|^2$ with a non-positive left-hand side and a non-negative right-hand side, resulting in $w^\top p\!=\!0$ and $Q p=\lambda_{\min}(Q)p$. This proves the sufficient part. Let $r\in\mathbb{R}^n$ be an eigenvector with respect to the minimal eigenvalue of $Q$. If $w^\top r=0$ holds, we get $r_\mathrm{n}^\top\big(Q+ww^\top\big)r_\mathrm{n}=\lambda_{\min}\left(Q\right)$ with $r_\mathrm{n}=r/\|r\|$ and thus $g(1,Q,w)\leq \lambda_{\min}(Q)$. Hence, $g(\cdot,Q,w)$ being non-decreasing and $g(0,Q,w)=\lambda_{\min}(Q)$ imply the necessary part.
\end{pf}


\begin{lem}\label{lemma:WellPosednessLipschitzianFunction}
    Let $\!\mathscr{h}:[0,\infty)\to\mathbb{R}$ be Lipschitz continuous and $\mathcal{N}\!=\!\mathrm{arg}\max_{\phi\geq 0}\big\{\!\mathscr{h}(\phi)\big\}$. If $\mathcal{N}\!\neq\! \emptyset$, then $\inf \mathcal{N}\!=\!\min \mathcal{N}\!$.\vspace{-0.75mm}
\end{lem}
\begin{pf}
    Let $ \mathcal{N}\neq\emptyset$, $\phi^\star=\mathrm{inf }\mathcal{N}$, $k=\mathrm{max}_{\phi\geq 0}\!\big\{\!\mathscr{h}(\phi)\big\}-\!\mathscr{h}(\phi^\star)$ and $L>0$ be a Lipschitz constant. Seeking a contradiction, we assume $\phi^\star\notin \mathcal{N}$. Thus, $k>0$. For all $\phi''>\phi^\star$, there exists $\phi'\in \mathcal{N}$ with $\phi'\leq \phi''$ since $\phi^\star$ is the greatest lower bound of $\mathcal{N}$. Setting $\phi''\!=\!\phi^\star+(2L)^{-1}k$, we get $\phi'-\phi^\star\leq (2L)^{-1}k$ and $\!\mathscr{h}(\phi')=\max_{\phi\geq 0}\!\big\{\!\mathscr{h}(\phi)\big\}$. Hence, we obtain $k=\!\mathscr{h}(\phi')-\!\mathscr{h}(\phi^\star)\leq L\,(\phi'-\phi^\star)\leq 2^{-1}k$ by means of the Lipschitz continuity, which contradicts $k>0$.\vspace{-0.75mm}
\end{pf}
\begin{lem}\label{lemma:BlockMatrices}
    Consider the block matrix $X= \begin{pmatrix}
            M&m\\m^\top & 0
        \end{pmatrix}$ with $M\in\mathbb{R}^{n\times n}$ and $m\in\mathbb{R}^n$. Then, $X$ is negative semi-definite, i.e. $X\!\leq\! 0$, if and only if $M\!\leq\! 0$ and $m\!=\!0$ hold.\vspace{-0.75mm}
\end{lem}
\begin{pf}
    Let $X\leq 0$. We get $v_1^\top M v_1+2v_1^\top m v_2\leq 0$ for all $v_1\in\mathbb{R}^n$, $v_2\in\mathbb{R}$. Now, $M\leq 0$ follows from $v_2=0$. Suppose, $m\neq 0$ holds. The choice $v_1=m$, $v_2=\|v_1\|^{-2}(1-v_1^\top M v_1)$ yields the contradiction $v_1^\top M v_1+2v_1^\top m v_2=1$. The other direction is trivial.\vspace{-0.75mm}
\end{pf}
\textbf{Matrices of the 
 Example}
\begin{align*}
    A=\begin{pmatrix}
        0&1\\-\omega_0^2 &-2\zeta\omega_0
    \end{pmatrix}\!, B=\begin{pmatrix}
        0 \\ \omega_0^{-2}
    \end{pmatrix}\!, P=\begin{pmatrix}3.9598 & 1\\
    1  &  1.4142\end{pmatrix}.
\end{align*}

\vspace{-1.25mm}
\bibliography{Bibliothek}



\end{document}